\documentclass[12pt]{article}
\usepackage{graphicx}
\usepackage{amsmath}
\usepackage{amssymb}
\usepackage{amsthm}

\newtheorem{teo}{Theorem}

\newtheorem{defi}{Definition}

\newtheorem{remark}{Remark}

\newcommand{\fin}{\end{document}}

\begin{document}

\title{\vspace{-1in}\parbox{\linewidth}{\footnotesize\noindent}
 \vspace{\bigskipamount} \\
Initial value problems in Clifford-type analysis
\thanks{ 
\hfil\break\indent
{\em $\dagger$E-mail: bolivarcolon@hotmail.com ~
\em $\ddagger$ corresponding author, e-mail: cvanegas@usb.ve }}}

\author{Yanett M. Bol\'{\i}var$\dagger$ \\
\emph{\small Universidad de Oriente, Departamento de Matem\'atica, Cuman\'{a}-Venezuela}\\
and \\
Carmen J. Vanegas$\ddagger$ \\
\emph{\small Universidad Sim\'on Bol\'{\i}var, Departamento de Matem\'atica,
Caracas-Venezuela}
}

\date{}

\maketitle

\begin{abstract} 

We consider an initial value problem of type 
$$
\frac{\partial u}{\partial t}={\cal F}(t,x,u,\partial_j u),
\quad u(0,x)=\varphi(x),
$$
where $t$ is the time, $x \in \mathbb{R}^n $ and $u_0$ is
a Clifford type algebra-valued function
satisfying
${\bf D}u=\displaystyle\sum_{j=0}^{n}\lambda_j(x)e_j\partial_ju = 0$,
$\lambda_j(x)\in \mathbb{R} $ for all $j$. We will solve this problem using 
the technique of associated spaces. In order to do that,
we give sufficient conditions on the coefficients of the operators
${\cal F}$ and ${\bf D}$, where  ${\cal F}(u)=
\displaystyle\sum_{i=0}^{n}A^{(i)}(x)\displaystyle\partial_iu$ for
$A^{(i)}(x) \in \mathbb{R}$ or $A^{(i)}(x)$ belonging to a Clifford-type
algebra, such that these operators are an associated pair.

\end{abstract}


\section{Introduction}
We consider the initial value problem
\begin{equation}\label{1ec1}
\frac{\partial u}{\partial t}={\cal F}(t,x,u,\partial_j u);
\end{equation}
\begin{equation}\label{1ec2}
u(0,x)=\varphi(x),
\end{equation}
where $t$ means the  variable time, $x=(x_0, \ldots, x_n)$ is a spacelike variable on ${\bf R}^{n+1}$
and $\partial_j$ is the operator diferentiation respect to $x_j$. 
The problem (\ref{1ec1}),(\ref{1ec2}) is equivalent to
the integro-differential
equation (see \cite{nagumo})
\begin{equation}\label{ec3}
u(t,x)=\varphi(x)+\int_0^t{\cal F}(\tau,x,u(\tau,x),\partial_j
u(\tau,x))d\tau.
\end{equation}
Therefore the solutions of the problem (\ref{1ec1}),(\ref{1ec2}) can be constructed
as fixed points of the operator
\begin{equation}\label{eco3}
U(t,x)=\varphi(x)+\int_0^t{\cal F}(\tau,x,u(\tau,x),\partial_j
u(\tau,x))d\tau.
\end{equation}

It is well known that the classical Cauchy-Kovalevskaya theorem gives a unique solution
to the problem (\ref{1ec1}),(\ref{1ec2}) considered in the context of complex analysis:
\begin{equation}\label{1ecc1}
\frac{\partial u}{\partial t}={\cal F}(t,z,u,\partial_z u);
\end{equation}
\begin{equation}\label{1ecc2}
u(0,z)=\varphi(z),
\end{equation}
provided ${\cal F}(t,z,u,\partial_z u)$ and $\varphi(z)$ are holomorphic functions in its variables.
The complex valued solution $u(t,z)$ of (\ref{1ecc1}), (\ref{1ecc2}) is holomorphic in $z$ and uniquely determined. 
H. Lewy showed that this problem has no solution if ${\cal F}$ is not a holomorphic function. He constructed
functions $f(t,x,y)$ infinitely many differentiable such that the equation:
$$
2i(x + iy)\partial_t w = \partial_x w + i\partial_y w + f(t, x, y)
$$
has no solution (see \cite{L}).
This shows that the equation (\ref{ec3}) does not always have a solution 
even if ${\cal F}(t,x,u,\partial_j u)$ and $\varphi(x)$ are
infinitely many differentiables. 

The concept of associated spaces \cite{heertush, tuts3, tuts5} leads to conditions under which the equation (\ref{ec3})
has solution. This concept comes from complex analysis:
In the holomorphic case (\ref{1ecc1}), (\ref{1ecc2}) the associated space is the space of holomorphic functions
and the right hand side ${\cal F}(t,z,u,\partial_z u)$ transforms this space into itself.

In order to apply a fixed-point theorem, the operator (\ref{eco3}) has to be estimated in a suitable function
space whose elements depend on $t$ and $x$. This can be done by
the so called interior estimate for the associated space.
Such estimate describes the behaviour of the derivatives near the boundary. In case of holomorphic functions,
such estimates can be obtained by the Cauchy Integral Formula.
Similar estimates can be shown in the framework of Clifford analysis.
This makes it possible to solve initial value problems with monogenic initial functions (see \cite{U}).
For initial value problems with other kind of initial functions see \cite{T, YC}.

Let ${\cal F}$ be a differential operator  while ${\cal G}$
is a differential operator with respect to the space variable $x$ whose coefficients do not depend on the time $t$.
The operators ${\cal F}$ and ${\cal G}$ are said to be associated  if ${\cal G}u=0$ implies
$ 
{\cal G}({\cal F}u)=0, ~ \mbox{for each} ~ t.
$
The  solutions of the differential equation ${\cal G}u=0$ form a function space called associated space to ${\cal G}$.
In \cite{{tuts3,tuts8,tuts5}} we can see that, if the 
initial function (\ref{1ec2}) satisfies an associated equation ${\cal G}u=0$
and the elements of the associated space satisfy an interior estimate,
then there exists a (uniquely determined) solution of the initial problem (\ref{1ec1}), (\ref{1ec2}) also 
satisfying the associated equation for
each $t$. 

In this paper, we will use this technique in order to solve the initial problem (\ref{1ec1}), (\ref{1ec2}) when
the operator ${\cal G}$ is given by
\begin{equation}\label{1ec5}
{\bf D}u=\sum_{j=0}^{n}\lambda_j(x)e_j\partial_ju.
\end{equation}
where the functions $\lambda_j$'s are supposed to be real valued
and $u$ is a  continuously differentiable function taking values in a Clifford type algebras depending on parameters. 
We call this operator the  generalized Cauchy-Riemann operator in
R$^{n+1}$.
First we will determine all ${\cal F}$ in the form
\begin{equation}\label{1par3}
{\cal F}u= \sum_{i=0}^{n}A^{(i)}(x)\partial_iu
\end{equation}
for which solutions $\varphi$ of $ {\bf D}{\cal \varphi}=0$ are admissible initial functions.
The functions $A^{(i)}(x)$ 
$i=0,1,\ldots n$, are functions of class ${\bf C}^1$ in a domain $\Omega$ of $ \mathbb{R}^{n+1}$
and are considered as real valued functions and as Clifford-type algebra valued.
Later on, the coefficients (belonging to $\mathbf{R}$) of the operator ${\cal F}$
are given and it is possible to obtain 
the conditions over the operator ${\bf D}$ such that the pair $({\cal F},{\bf D})$ is associated.


\section{Preliminaries}

A Clifford algebra depending on parameters (see \cite{tuts7}) can be defined as equivalence classes in the ring $R[X_{1},...,X_{n}]$ 
of polynomials in $n$ variables $X_{1},...,X_{n}$ with real coefficients, where two polynomials are said to be equivalent if 
their difference is a polynomial for which each term contains at least one of the factors
\begin{equation}\label{cesar-A}
X_{j}^{k_{j}}+\alpha_{j} ~ \mbox{and} ~ X_{i}X_{j}+X_{j}X_{i}-2\gamma_{ij},
\end{equation}
where $i,j=1, \cdots, n$ , $ i\neq j$, and
the $k_{j}\ge 2$ are natural numbers. 
The parameters $\alpha_{j}$ and $\gamma_{ij}=\gamma_{ji}$ have to be real and may depend also on further variables 
such as the variable $x$ in $\mathbf{R}^{n+1}$. If the parameters do not depend on further variables and if $n\geq 3$, 
the Clifford type algebra generated by the structure polynomials (\ref{cesar-A}) is denoted 
by ${\cal A}_{n}(k_{j},\alpha_{j},\gamma_{ij})$. For $n=1$ we write ${\cal A}_{1}(k,\alpha)$. 
In ${\cal A}_{n}(2,\alpha_{j},\gamma_{ij})$ all of the $k_{j}$ are equal to $2$ and 
$\alpha_{j}$ and $\gamma_{ij}$ are constant. This algebra has the dimension $2^n$.
The classical Clifford algebra 
${\cal A}_n(2, 1, 0)$ is denoted by ${\cal A}_n$. 

In this paper the functions $u$ will be consider in a domain $\Omega$ of  ${\bf R}^{n+1}$ with va\-lues in  
the particular algebra ${\cal
A}_n(2,\alpha_j,\gamma_{ij})$, thus $u(x)=\displaystyle\sum_{A\in \Gamma}u_A(x)e_A,$ 
where $\Gamma=\{0,1,\ldots,n,12,\ldots,12\ldots n\}$ and
each $u_A$ is a real valued function.

We will understand that a 
function $u(x)=\displaystyle\sum_{A\in \Gamma}u_A(x)e_A,$ with values in an
algebra de Clifford ${\cal A}_n(2,\alpha_j,\gamma_{ij})$ is of class {\bf C}$^1$ (or {\bf C}$^2$)
in the domain $\Omega$ if each $u_A$ is  of class {\bf C}$^1$ (or {\bf C}$^2$) in $\Omega$.

In similar form as the definition of monogenic functions  in the
classical algebra ${\cal A}_n$ (see \cite{Wolyu1}), we have
\begin{defi} A continuously differentiable function $u(x)$ with values in the generalized
Clifford algebra given by (\ref{cesar-A}) satisfying $ {\bf D}u=0$  is called generalized (left) monogenic function.
\end{defi}

For the functions $u$ of class {\bf C}$^1$ in $\Omega$, the conjugated
operator associated to generalized Cauchy-Riemman operator (\ref{1ec5})  is defined by
$$
\overline{\bf D}u=\lambda_0(x)\partial_0u-\sum_{j=1}^{n}\lambda_j(x)e_j\partial_ju.
$$
Since $\lambda_i(x)$'s are supposed to be real valued, a monogenic function $u$ of class {\bf C}$^2$ verifies
\begin{eqnarray}\label{acopla}
{\bf \overline{D}D}u
&=&\sum_{i=0}^{n}\lambda_0(x)\partial_0(\lambda_i(x))e_i\partial_iu -
\sum_{j=1}^n\sum_{i=0}^{n}\lambda_j(x)e_j\partial_j(\lambda_i(x))e_i\partial_iu \nonumber \\
 &+&\lambda_0^2(x)\partial_0^2u+\sum_{i=0}^{n}\alpha_i \lambda_i^2(x) \partial_i^2u
-2\sum_{i<j}\gamma_{ij}\lambda_i(x)\lambda_j(x)\partial_i\partial_ju=0.
\end{eqnarray}
If $\lambda_i(x)$ are real constants for all $i= 1, \ldots, n$, then 
the partial differential equation (\ref{acopla}) is elliptical provided 
\begin{equation}\label{C1}
\alpha_j>0 ~~ \mbox{and} ~~ |\gamma_{ij}\lambda_i\lambda_j|\leq k, ~~  i,j=1,\ldots n,
\end{equation}
for a suitable constant $k$. Further if the underlying algebra is ${\cal A}_n(2, \alpha_j, 0)$
and $\lambda_i(x)$'s are real valued functions, then (\ref{acopla}) is also an elliptic equation  under the condition
\begin{equation}\label{C2}
\alpha_j>0 ~~ j=1,\ldots n.
\end{equation}

\subsection {A formula for the product {\bf D}$(u.v)$}

Now, based on the previously seen definitions, a formula
of the generalized Cauchy-Riemman operator applied to the product
of two functions $u$  and $v$ with values in ${\cal A}_n(2,\alpha_j,\gamma_{ij})$ is built.

Let $u$ and $v$ be functions defined in $\mathbb{R}^{n+1}$, continuously diferentiable and with values in ${\cal A}_n(2,\alpha_j,\gamma_{ij})$. 
Using induction over $n$, we shall show the identity
\begin{equation}\label{lashpro}
\partial_k(u\cdot v)=\partial_k(u)\cdot v+u\cdot \partial_k(v),\; k=0,\ldots,n.
\end{equation}

For the case $n=1$. Let $\Gamma_0=\{0,1\}$ be and $u=\displaystyle\sum_{i\in\Gamma_0}u_i e_i$,
$v=\displaystyle\sum_{i\in\Gamma_0}v_i e_i$. Then $u\cdot v = \displaystyle\sum_{i,j\in\Gamma_0}u_iv_je_{ij}$
and trivially $ \partial_k(u\cdot v) = \partial_k(u)\cdot v+u\cdot \partial_k(v)$.

Now we suppose that the identity (\ref{lashpro})  is valid for $n=k$.
Let's consider the index sets 
$\Gamma=\{0,1,2,\ldots,12,13,\ldots,123\ldots k+1\}\,,$  

\noindent $\Gamma_1=\{0,1,2,\ldots,12,13,\ldots,123\ldots k\}$  and 
$\Gamma_2=\Gamma-\Gamma_1.$ 
Thus
 $u=u^{(1)}+u^{(2)}$ and $v=v^{(1)}+v^{(2)}$
where $u^{(i)},\;v^{(i)}$ are a linear combination of the elements $e_A$ of 
the bases indexed by the sets $\Gamma_i$, $i=1,2.$
Then
\begin{eqnarray*}
\partial_k(u\cdot v)   &=&\partial_k[\left(u^{(1)}+u^{(2)}\displaystyle \right)(v^{(1)}+v^{(2)})] \\
&=&\partial_k(u^{(1)}v^{(1)})+\partial_k(u^{(1)}v^{(2)})
+\partial_k(u^{(2)}v^{(1)})+\partial_k(u^{(2)}v^{(2)}).
\end{eqnarray*}
Applying the induction hypothesis and rearrange the terms we obtain that
 $\partial_k(u\cdot v) =\partial_k(u)\cdot v+u\cdot \partial_k(v)$
and the proof is obtained.
\\
Now, according to definition (\ref{1ec5}),
\begin{eqnarray}\label{1pro3}
{\bf D}(u\cdot v)&=&\sum_{i=0}^{n} \lambda_i(x)e_i\partial_i(u\cdot v)
=\sum_{i=0}^{n} \lambda_i(x) e_i(\partial_i(u)\cdot v+u\cdot \partial_i(v)) \nonumber \\  
&=& {\bf D}(u)\cdot v +\sum_{i=0}^{n}\lambda_i(x)e_i(u\cdot \partial_i(v)).
\end{eqnarray}
If $n=2$, the general formula (\ref{1pro3}) yields
\begin{eqnarray}\label{1pro2,2}  
{\bf D}(u\cdot v)&=& {\bf D}(u)\cdot v+u\cdot{\bf D}v  \\
\nonumber &+&2\lambda_1[- u_2 \gamma e_0 - u_{12} \gamma e_1 - u_{12} \alpha_1e_2+ u_2e_{12}]\partial_1v\\
\nonumber  &+&2\lambda_2[ u_1 \gamma e_0+  u_{12} \alpha_2e_1+
u_{12} \gamma e_2 - u_1e_{12}]\partial_2v.
\end{eqnarray}
This expression is a consequence of using
$$
e_1\cdot u=u\cdot e_1+2[- u_2 \gamma e_0 - u_{12} \gamma e_1 - u_{12} \alpha_1e_2+ u_2e_{12}]
$$
and
$$
e_2\cdot u=u\cdot e_2+2[ u_1e_0 \gamma+ u_{12} \alpha_2e_1+ u_{12} \gamma e_2 - u_1e_{12}].
$$

\section{Sufficient conditions}
\subsection{Conditions over the coefficients of ${\cal F}$}
We consider the operator  ${\cal F}(u)$ defined by (\ref{1par3}), where
$u$ is a ${\cal A}_n(2,\alpha_j,\gamma_{ij})$-valued function and
we will determine conditions over  $A^{(i)}$  guaranteeing
\begin{equation}\label{1ec4} 
{\bf D}u=0 \Rightarrow {\bf D}({\cal F}u)=0.
\end{equation}
From the equation ${\bf D}u=0$  we obtain
$$
\partial_0u(x) = -\sum_{j=1}^{n}\beta_j(x)e_j\partial_ju,
$$
where $\beta_j(x)=\displaystyle\frac{\lambda_j(x)}{\lambda_0(x)},\;j=1,\ldots,n$,
$\lambda_0(x)\neq 0$. This formula leads to the equality
$$
\partial_k\partial_0u(x) =-\sum_{j=1}^{n}\partial_k(\beta_j(x)\cdot e_j\cdot \partial_ju(x))\quad k=0,1,\ldots n.
$$
From $u\in C^2$, we get $\partial_k\partial_0u=\partial_0\partial_ku, \; k=1,\ldots n.$
 Thus  both last equations make this way possible to write
\begin{equation}\label{3ge-mono} 
{\bf D}(\partial_ku(x)) =-\sum_{j=1}^{n}\lambda_0\partial_k(\beta_j(x)).e_j.\partial_ju(x), \quad  k=0,1, \ldots  n.
\end{equation}
The different conditions to be set will depend on the Cauchy-Riemman operator chosen and of the
characteristics of $A^{(i)}(x)$'s.

{\bf Case I: $A^{(i)}$'s are real valued functions}.~
Applying $D$ to (\ref{1par3}) and considering (\ref{3ge-mono}), it follows that 
${\bf D}({\cal F}u)$ can be expressed as a linear combination of the first order derivatives of $u$:
\begin{eqnarray*}
{\bf D}({\cal F}u) &=&     
\sum_{i=0}^{n}\left({\bf D}(A^{(i)}(x))\cdot \partial_iu + A^{(i)}(x)\cdot {\bf D}(\partial_iu )\right) \\ \nonumber 
&=&\sum_{i=1}^{n}({\bf D}A^{(i)}-{\bf D}A^{(0)}\beta_ie_i)\displaystyle\partial_i u 
-\sum_{i=1}^{n}\sum_{j=0}^{n}A^{(j)}(x)\lambda_0 \partial_j(\beta_i) e_i \partial_iu
\end{eqnarray*}
and equating its coefficients to zero we obtain the following $n$ sufficient conditions over the $n+1$
real functions $A^{(i)}$:
\begin{equation}\label{consufi1}
{\bf D}A^{(i)}(x)-{\bf D}A^{(0)}(x)\beta_i(x) e_i -\sum_{j=0}^{n}A^{(j)}(x)\lambda_0(x) \partial_j\beta_i (x) e_i = 0,
\end{equation}
for $i = 1, \cdots, n$. These conditions guarantee (\ref{1ec4}).

{\bf Case II: $A^{(i)}$'s are ${\cal A}_n(2,\alpha_j,\gamma_{ij})$-valued functions}.~
Suppose  $A^{(i)}$ functions of class {\bf $C^1$}. 
Due to the heavy calculations to be made, we will restrict to the case $n = 2$. 
 From (\ref{1par3}) and $n=2$ we have
$ {\cal F}u= \sum_{i=0}^{2}A^{(i)}(x)\partial_iu,$
and since $A^{(i)}$'s  are defined in  ${\bf R^3}$ with values in ${\cal A}_2(2,\alpha_1,\alpha_2,\gamma)$ 
formula (\ref{1pro2,2}) is used

\begin{eqnarray*}
\label{2pro-2} {\bf D}({\cal F}u)
 &=&\sum_{i=0}^{2}{\bf D}\left[A^{(i)}(x).\partial_iu  \right]\\
\nonumber   &=& \sum_{i=0}^{2}{\bf
D}\left[A^{(i)}(x)\right].\partial_iu +A^{(i)}(x).{\bf D}
\left[\partial_iu  \right]\\
\nonumber   &+&2\lambda_1\left[- A^{(0)}_2 \gamma e_0 - A^{(0)}_{12}
\gamma e_1 - A^{(0)}_{12} \alpha_1e_2+  A^{(0)}_2e_{12}
\right]\partial_1\partial_0 u\\
\nonumber & +&2\lambda_2\left[  A^{(0)}_1 \gamma e_0+  A^{(0)}_{12}
\alpha_2e_1+ A^{(0)}_{12} \gamma e_2 -
A^{(0)}_1e_{12}\right]\partial_2\partial_0 u.
  \\
\nonumber   &+&2\sum_{i=1}^{2}\left\{\lambda_1\left[- A^{(i)}_2
\gamma e_0 - A^{(i)}_{12} \gamma e_1 - A^{(i)}_{12} \alpha_1e_2+
A^{(i)}_2e_{12}\right]\right.\partial_1\partial_i
u  \\
\nonumber &+&\lambda_2\left.\left[  A^{(i)}_1 \gamma e_0+
A^{(i)}_{12} \alpha_2e_1+  A^{(i)}_{12} \gamma e_2 -
A^{(i)}_1e_{12}\right]\partial_2\partial_i u\right\},
\end{eqnarray*}
where $A^{(k)}= A_0^{(k)} + A_1^{(k)}e_1 + A_2^{(k)}e_2 +  A_{12}^{(k)}e_{12}$.
From (\ref{3ge-mono})
\begin{eqnarray}
\label{pro} {\bf D}({\cal F}u) \nonumber &=&\sum_{i=1}^{2}({\bf D}A^{(i)}-{\bf D}A^{(0)}\beta_ie_i)\displaystyle\partial_i u
\\
 \nonumber &-&\sum_{i=1}^{2}\sum_{j=0}^{2}A^{(j)}(x)\lambda_0.
 \partial_j(\beta_i).e_i.\partial_iu
\\
 \nonumber &-&2\lambda_1\left[- A^{(0)}_2
\gamma e_0 - A^{(0)}_{12} \gamma e_1 - A^{(0)}_{12} \alpha_1e_2+
A^{(0)}_2e_{12}
\right]\\
 \nonumber&&[\sum_{i=1}^{2}\partial_1(\beta_i).e_i.\partial_iu+\beta_i.e_i.\partial_1(\partial_iu)]\\
\nonumber & -&2\lambda_2\left[  A^{(0)}_1 \gamma e_0+  A^{(0)}_{12} \alpha_2e_1+ A^{(0)}_{12} \gamma e_2 -
A^{(0)}_1e_{12}\right]\\
 \nonumber&&[\sum_{i=1}^{2}\partial_2(\beta_i).e_i.\partial_iu+\beta_i.e_i.\partial_2(\partial_iu))].
  \\
\nonumber   &+&2\sum_{i=1}^{2}\left\{\lambda_1\left[- A^{(i)}_2
\gamma e_0 - A^{(i)}_{12} \gamma e_1 - A^{(i)}_{12} \alpha_1e_2+
A^{(i)}_2e_{12}\right]\right.\partial_1\partial_i
u  \\
\nonumber &+&\left.\lambda_2\left[  A^{(i)}_1 \gamma e_0+
A^{(i)}_{12} \alpha_2e_1+  A^{(i)}_{12} \gamma e_2 -
A^{(i)}_1e_{12}\right]\partial_2\partial_i u\right\}.
\end{eqnarray}
 If  $u$  is assumed of class $C^2$ then the condition (\ref{1ec4}) is verified provided the
coefficients of the first order derivatives satisfy
\begin{eqnarray}
\label{1no}({\bf D}A^{(1)}-{\bf D}A^{(0)}\beta_1e_1) -\sum_{i=0}^{2}A^{(i)}(x)\lambda_0
 \partial_i\beta_1.e_1=0,   \\
 \label{1no1} ({\bf D}A^{(2)}-{\bf D}A^{(0)}\beta_2e_2)
-\sum_{i=0}^{2}A^{(i)}(x)\lambda_0
 \partial_i\beta_2.e_2=0,
\end{eqnarray} 
producing altogether 8 equations. Additional conditions for second order
coe\-ffi\-cients are
\begin{equation}
\begin{array}{lll}\label{1no2}
  \alpha_1\beta_1A^{(0)}_{12}  = A^{(1)}_2&& -  \beta_1A^{(0)}_2 =
     A^{(1)}_{12}, \\\\
      \alpha_2\beta_2 A^{(0)}_{12}=- A^{(2)}_1  &&
 \beta_2 A^{(0)}_1 =
     A^{(2)}_{12},\\\\
      2 \beta_1\beta_2\gamma A^{(0)}_{12}=\beta_2A^{(1)}_1  - \beta_1 A^{(2)}_2,
\end{array}
\end{equation}
 and they allow reducing the number of
unknown quantities in (\ref{1no}) and (\ref{1no1}) to 21.
 Finally, conditions (\ref{1no})-(\ref{1no2}) guaranteeing that the pair of operators
${\cal F}u = \sum_{i=0}^{2}A^{(i)}(x)\displaystyle\partial_iu$ and
${\bf D}u = \sum_{j=0}^{2}\lambda_j(x)e_j\partial_ju$ be associated,
have been obtained.
 For the classic algebra ${\cal A}_2$ and the case in which all $\lambda_i$ are equal to $1$ and the components
$A^{(i)}_{j}$, $i=1,\ldots, n$, $j=1,\ldots, n,$ linear functions depending of the variables $x_0$, $x_1$ and $x_2$, of
(\ref{1no}) and (\ref{1no1}), 21-8=13 linearly independent associated operators ${\cal F}$, 
called admissible, are obtained. 
\begin{teo}\label{coefF} Consider the operator ${\bf D}$ defined by (\ref{1ec5})
where the coefficients $\lambda^{(i)}$ are real valued functions and of class ${\bf C}^1$
defined in $\Omega \subset \mathbb{R}^{n+1}$. Then ${\bf D}$  
is associated to each operator ${\cal F}$ given by (\ref{1par3})
if the conditions (\ref{consufi1}) or (\ref{1no})-(\ref{1no2}) are satisfied. 
\end{teo}
\begin{remark}
The results given in Theorem (\ref{coefF}) are an extension of those shown 
in \cite{2sontush} for the classic algebra ${\cal A}_2$.
\end{remark}

\subsection{Conditions over the coefficients of  ${\bf D}$}

Now we consider the operator  ${\cal F}$ given by (\ref{1par3}) and suppose the coefficients $A^{(i)}$ are real valued.
In order to determine conditions over the real valued functions $\lambda_i$
such that (\ref{1ec4}) be valid, we asume $u$ as a monogenic function with values in ${\cal A}_{n}(2,\alpha_j,\gamma_{ij})$.
Then we obtain
\begin{equation*} \label{pro6} 
{\bf D}({\cal F}u) = \sum_{i=1}^{n}({\bf D}A^{(i)}-{\bf D}A^{(0)}\cdot \beta_ie_i)\partial_i u -
\sum_{i=1}^{n}\sum_{j=0}^{n}A^{(j)}(x)\lambda_0 \partial_j(\beta_i)e_i \partial_iu. 
\end{equation*}
Therefore the operator ${\cal F}$ is associated to operator $D$ if 
the coefficient of each  $\partial_i u$ vanish identically, i.e., 
$$
({\bf D}A^{(i)}-{\bf D}A^{(0)}\cdot \beta_ie_i) - 
\sum_{j=0}^{n}A^{(j)}(x)\lambda_0 \partial_j(\beta_i) e_i = 0 ~~ \mbox{for each} ~~ i= 1, \cdots, n.
$$
Dividing by $\lambda_0$ and applying operator ${\bf D}$ we have the system
\begin{eqnarray}
&& \partial_0 A^{(i)} + \alpha_i \beta_i^2 \partial_i A^{(0)} - 2 \beta_i \sum_{k=i+1}^n {\gamma_{ik}\beta_k\partial_k A^{(0)}} = 0 \label{betai1}\\
&& \beta_j \partial_j A^{(i)} = 0, ~~ \mbox{for} ~ j \ne i \label{betai2}\\
&& \beta_i \partial_i A^{(i)} - \beta_i \partial_0 A^{(0)} - \sum_{j=0}^n {A^{(j)}\partial_j(\beta_i)} = 0 \label{betai3}\\
&& \beta_j \beta_i \partial_j A^{(0)} = 0, ~~ \mbox{for} ~  j \ne i \label{betai4}.
\end{eqnarray}
If $\beta_j, \beta_i\neq 0$, the system  (\ref{betai1})-(\ref{betai4}) leads to an uncoupled system of equations
\begin{equation}\label{conA1} 
\beta_i(\partial_i A^{(i)} - \partial_0 A^{(0)}) - \sum_{j=0}^n {A^{(j)}\partial_j(\beta_i)} = 0.
\end{equation}
Therefore we have proved the following
\begin{teo}\label{coefA} Consider the operator ${\cal F}$ defined by (\ref{1par3})
with real valued coefficients $A^{(i)}= A^{(i)}(x_i)$ arbitrarily given and 
of class ${\bf C}^1$, then ${\cal F}$  
is associated to each generalized Cauchy- Riemman operator given by (\ref{1ec5})
for which the functions $\beta_i=\displaystyle\frac{\lambda_i}{\lambda_0}$, $\lambda_0\neq 0$
satisfy the system (\ref{conA1}).
\end{teo}
\begin{remark} The former theorem generalizes the results appearing in Theorem 3 of  \cite{tuts8}
in which the author tried the case $n= 3$ in ${\cal A}_n$.
\end{remark}

\section{Initial value problems}

Using the theory of associated spaces  \cite{heertush, tuts3, tuts5} the initial value problem
(\ref{1ec1}), (\ref{1ec2})
can be solved.
In this work, conditions for the associated pair $({\cal F},{\bf D})$ have been obtained in two ways:
Given ${\bf D}$, if the equations (\ref{consufi1}) or (\ref{1no})-(\ref{1no2}) are satisfied, 
we can get the conditions
on the functions $A^{(i)}$ which allow to find
${\cal F}$, reciprocally given ${\cal F}$, through (\ref{conA1}),
we can determine the operator  ${\bf D}$.
Thus one obtains immediately that the operator
${\cal F}$ sends monogenic functions into monogenic functions.

It is easy to see that solutions of the initial value problem (\ref{1ec1}), (\ref{1ec2}) are fixed points 
of the integro-diferential operator (\ref{eco3})
and vice versa.
In order to apply a fixed point theorem, solutions of ${\bf D}u = 0 $ must satisfy
an interior estimate of first order (see \cite{2sontush, tuts8, tuts5}).
Therefore the  first order
derivatives of the solution $u$, contained in the operator (\ref{eco3})   
can be estimated.
This requirement is achieved if the solutions of ${\bf D}u=0$
are solutions of an elliptical differential equation (see \cite{niremberg, tuts5}).
In our case, such solutions satisfy the equation (\ref{acopla}) which is 
an elliptic equation al least under the condition (\ref{C1}) or (\ref{C2}). 
Hence we can apply the following theorem (see \cite{tuts3, tuts5})
\begin{teo}  Suppose ${\cal F}$ and ${\bf D}$ are a pair of associated operators 
and the solutions of the associated equation ${\bf D}u=0$ satisfy an interior estimate of first order.
Then the initial value problem $(\ref{1ec1})-(\ref{1ec2})$ is soluble provided 
the initial solution $\varphi$ satisfies the condition ${\bf D}\varphi=0$.
\end{teo}


\begin{thebibliography}{99}


\bibitem{niremberg} Agmon S., Douglis A., Niremberg L., Estimates near the boundary for solutions of elliptic partial differential
equations satisfying general boundary conditions I and II.  {\it Comm. Pure Appl. Math.,} vol. 12, pp. 623-727, 1959, and vol. 17,
pp. 35-92, (1964).

\bibitem{heertush}  Heersink R., Tutschke W.,
Solution of initial value problems in associated spaces.
{\it Functional Analytic Methods in Complex Analysis and Applications to Partial
Di\-ffe\-ren\-tial Equations}.  World Sci. Publ. pp. 209-219, (1995).

\bibitem{L}
Lewy H., An example of a smooth linear partial differential
equation without solution. {\it Ann. of Math.}, vol. 66, pp. 155-158, (1957).

\bibitem{nagumo} Nagumo M., \"{U}ber das Anfangswertproblem Partieller  
Di\-ffe\-ren\-tial\-glei\-chun\-gen, {\it Japan.J.Math.}, vol.
18, pp. 41-47, (1941-43).

\bibitem{2sontush} Son L., Tutschke W.,
Complex methods in higher dimensions- recent trends for solving boundary
value and initial value problems. {\it Complex Variables}, vol. 50, No. 7-11, pp. 673-679, (2005).

\bibitem{T} Tutschke W.,  Solution of initial value problem in classes of generalized analytic functions.
Teubner Leipzig and Springer Verlag, (1989).

\bibitem{tuts3}   ------, 
The method of weighted function spaces for solving initial value and boundary value problems.
Contained in {\it Functional-analytic and complex methods, their interactions and
aplications to partial differential equations}.  World Scientific,  pp. 75-90, (2001).

\bibitem{tuts8}     ------,
Complex analysis within the framework of analytical methods for partial differential equations. 
Contained in {\it Methods of Complex and Clifford Analysis}. Proceedings of the International Conference on Applied
Mathematics based on partial differential equations and complex analysis (ICAM), Hanoi August 25-29, 2004. SAS International
Publications, pp. 123-140, (2006).
  
\bibitem{tuts5}  ------,
Associated spaces - a  new tool of real and complex analysis. {\it Natl. Univ. Publ. Hanoi}, (2008).

\bibitem{Wolyu1} Tutschke W., Vanegas C.,
 M\'etodos del an\'alisis complejo en dimensiones superiores,
{\it XXI Escuela Venezolana de Matem\'aticas}. Ediciones IVIC, (2008).

\bibitem{tuts7} ------,
Clifford algebras depending on parameters and their applications to partial differential equations. Contained in  
{\it Some topics on value distribution and differentiability in complex and p-adic analysis}.
Science Press Beijing, pp. 430-450, (2008).

\bibitem{U}  Y\"{u}ksel U., Solution of initial value problems with monogenic initial functions
in Banach Spaces with $L_p$-norm. {\it Adv. appl. Clifford alg.}, vol.20, No.1, pp. 201–209, (2010).

\bibitem{YC} Y\"{u}ksel, Celebi O., Solution of Initial Value Problems of Cauchy-Kovalevsky Type in
the Space of Generalized Monogenic Functions. {\it Adv. appl. Clifford alg}, vol. 20, No.2, pp. 427-444, (2010).

\end{thebibliography}
\end{document}